\documentclass[12pt]{article}

\usepackage{amssymb}

\newtheorem {problem}{\bf Problem}
\newtheorem {definition}{\bf Definition}

\newtheorem {lemma} {\bf Lemma}
\newtheorem{corollary} {\bf Corollary}
\newtheorem {theorem} {\bf Theorem}

\title{Fibonacci sequence related to a combinatorial problem on binary matrices}
\author{Krasimir Yordzhev}
\date{}

\begin{document}
\maketitle

\begin{center}Faculty of Mathematics and Natural Sciences, South-West University

2700 Blagoevgrad, Bulgaria

E-mail: yordzhev@yahoo.com
\end{center}

\begin{abstract}
We discuss an equivalence relation on the set of square binary matrices with the same number of 1's in each row and each column. Each binary matrix is represented using ordered $n$-tuples of natural numbers. We give a few starting values of integer sequences related to the discussed problem. The obtained sequences are new and they are not described in the On-Line Encyclopedia of Integer Sequences  (OEIS). We show a relationship between some particular values of the parameters and the Fibonacci sequence.
\end{abstract}

\textbf{Keyword:} binary matrix; equivalence relation; factor-set; Fibonacci number

\textbf{2010 Mathematics Subject Classification:}  05B20; 11B39

\section{Introduction}
A \emph{binary} (or   \emph{boolean}, or (0,1)-\emph{matrix}) is a matrix whose all elements belong to the set $\mathcal{B} =\{ 0,1 \}$. With $\mathcal{B}_n$ we will denote the set of all  $n\times n$  binary matrices.

Let $n$ and $k$ be positive integers. We let $\Lambda_n^k$ denote the set of all $n\times n$ binary matrices in each row and each column of which there are exactly $k$ in number 1's.
Let us denote with $\lambda (n,k) =|\Lambda_n^k |$ the number of all elements of $\Lambda_n^k$.
There is not any known formula to calculate  the $\lambda (n,k)$  for all $n$ and $k$.

Let $A,B\in\Lambda_n^k$. We will say that $A\sim B $, if $A$ is obtained from $B$ by moving some rows and/or columns. Obviously,  the relation  defined like that is an equivalence relation. We denote with
\begin{equation}\label{mu(n,k)}
\mu (n,k) = \left| {\Lambda_n^k}_{/_\sim} \right|
\end{equation}
the number of equivalence classes on the above defined relation.

\begin{problem}\label{prbl}
Find $\mu(n, k)$ for given integers $n$ and $k$, $1\le k<n$.
\end{problem}

Problem \ref{prbl} is the subject of discussion in this article.

\section{Some values of the integer function $\mu (n,k)$}

The task of finding the number of equivalence classes for all integers $n$ and $k$, $1\le k\le n$ is an open scientific problem. We partially solve this problem by implementing a computer program to find $\mu(n, k)$ for some (not great) values of $n$ and $k$. Moreover, using bitwise operations, our algorithm received one representative from each equivalence class without examining the whole set
$\Lambda_n^k$ \cite{umb2013}.

Let $\mathbb{N} $ be the set of natural numbers and let
\begin{equation}\label{T_n}
\mathcal{T}_n =\left\{ \langle x_1 ,x_2 ,\ldots ,x_n \rangle \; |\; x_i \in \mathbb{N},\; 0\le x_i \le 2^n -1,\; i=1,2,\ldots ,n\right\} .
\end{equation}

In \cite{umb2009} and \cite{Kostadinova} we describe an one-to-one correspondence
\begin{equation}\label{varphi}
\varphi \; :\; \mathcal{B}_n \cong \mathcal{T}_n
\end{equation}
which is based on the binary presentation of the natural numbers. If $A\in \mathcal{B}_n$ and $\varphi (A) =\langle x_1 ,x_2 ,\ldots x_n \rangle $, then $i$-th row of  $A$ is integer $x_i$ written in binary notation.

In \cite{Kostadinova} we prove that the representation of the elements of $\mathcal{B}_n$ using ordered $n$-tuples of natural numbers leads to making a fast and saving memory algorithms.

Let $A\in \mathcal{B}_n$ and let $\mathbf{x} =\langle x_1 ,x_2 ,\ldots ,x_n \rangle =\varphi (A)$. Then we denote
$$\mathbf{x}^t =\varphi (A^t ) ,$$
where $A^t \in \mathcal{B}_n$ is the transpose of the matrix  $A$ .

\begin{definition}\label{defcan}
Let $\mathbf{x}= \langle x_1 ,x_2 , \ldots ,x_n \rangle\in \mathcal{T}_n$  and let $\mathbf{x}^t = \langle y_1 ,y_2 ,\ldots ,y_n \rangle $. The element $\bf x\in \mathcal{T}_n$ we will call  {\bf canonical element}, if $x_1\le x_2 \le \cdots \le x_n$  and $y_1\le y_2 \le \cdots \le y_n$. The matrix $A\in \Lambda_m^k$  we will call  {\bf canonical matrix}, if $\varphi (A)$ is a canonical element in $\mathcal{T}_n$, where $\varphi$ is the defined with (\ref{varphi}) isomorphism.
\end{definition}

Obviously, when $k = 0$, the zero $n\times n$ matrix is the only matrix in the set $\Lambda_n^0$. When $k = n$, there is only one $n\times n$ binary matrix of $\Lambda_n^n$, and this is the matrix all elements of which are equal to 1. Therefore

\begin{equation}\label{keq0}
\mu (n,0)=\mu (n,n) =1 .
\end{equation}

It is easy to prove that for any positive integer $n$ is satisfied

\begin{equation}\label{keq1}
\mu (n,1)=\mu (n,n-1) =1 .
\end{equation}

When $k = 1$ the only canonical element is $\mathbf{x} = \langle 1,2,4, \ldots , 2^{n-1} \rangle \in \mathcal{T}_n$, i.e., if $A\in \Lambda_{n}^{1}$ is a canonical matrix, then $ A $ is a binary matrix with 1 in the second (not leading) diagonal and 0 elsewhere. For $k=n-1$, if $A\in \Lambda_{n}^{n-1}$ is a canonical matrix, then $A$ is a binary matrix with 0 in the leading diagonal and 1 elsewhere.

An algorithm for finding all canonical elements of $\mathcal{T}_n$ is described in detail in \cite{umb2013}. For $k = 2,3,4$ and $k = 5$, we will display the first elements of the sequences $\displaystyle \left\{ \mu (n,k) \right\}_{n=k}^\infty $ for some values of the parameter $n$. Using a computer program \cite{umb2013} we obtained the following  sequences
\begin{equation}\label{k2}
\left\{ \mu (n,2) \right\}_{n=2}^\infty \; =\; \left\{ 1,\; 1,\; 2,\; 5,\; 13,\; 42,\;155,\; 636,\; 2\; 889 ,\; 14\; 321 ,\; \ldots \right\}
\end{equation}
\begin{equation}\label{k3}
\left\{ \mu (n,3) \right\}_{n=3}^\infty \; =\; \left\{ 1,\; 1,\; 3,\; 25,\; 272,\; 4\; 070,\; 79\; 221, \; \ldots \right\}
\end{equation}
\begin{equation}\label{k4}
\left\{ \mu (n,4) \right\}_{n=4}^\infty \; =\; \left\{ 1,\; 1,\; 5,\; 161,\; 7\; 776  ,\; 626\; 649 ,\; \ldots \right\}
\end{equation}
\begin{equation}\label{k5}
\left\{ \mu (n,5) \right\}_{n=5}^\infty \; =\; \left\{ 1,\; 1,\; 8,\; 1\; 112  ,\; 287\; 311 ,\; \ldots \right\}
\end{equation}

The obtained integer sequences (\ref{k2}) $\div$ (\ref{k5}) are not described in the On-Line Encyclopedia of Integer Sequences  (OEIS)\cite{OEIS}.

\section{The function $\mu(n,k)$ and Fibonacci numbers}
The sequence $\displaystyle \left\{ f_n \right\}_{n=0}^\infty $ of Fibonacci numbers is defined by the recurrence relation (see for example \cite{atanassov} or \cite{koshy})

\begin{equation}\label{Fib}
f_0 =f_1 =1, \qquad f_n = f_{n-1} +f_{n-2} \quad \textrm{for} \quad n=2,3,\ldots
\end{equation}

In this section, we will prove that the sequence $\displaystyle \left\{ \mu (k+2,k) \right\}_{k=0}^\infty$ coincides with the Fibonacci sequence (\ref{Fib}).

\begin{lemma}\label{lemm1}
If $A=(\alpha_{i\, j} )\in \Lambda_n^k$ is a canonical matrix then
$$\alpha_{1\, 1} =\alpha_{1\, 2} =\cdots =\alpha_{1\; n-k} =0,\quad \alpha_{1\; n-k+1} =\alpha_{1\; n-k+2} =\cdots =\alpha_{1\, n} =1,$$
$$\alpha_{1\, 1} =\alpha_{2\, 1} =\cdots =\alpha_{n-k\; 1} =0,\quad \alpha_{n-k+1\; 1} =\alpha_{n-k+2\; 1} =\cdots =\alpha_{n\, 1} =1,$$
\end{lemma}

Proof. Immediately.

\hfill $\Box$

\begin{corollary}
If $\mathbf{x} =\langle x_1 ,x_2 ,\ldots ,x_n \rangle \in \mathcal{T}_n$ is a canonical element then $x_1 = 2^k -1$.

\hfill $\Box$
\end{corollary}

\begin{theorem}
Let the sequence $\displaystyle \left\{ \mu (k+2,k)\right\}_{k=0}^\infty$ is defined by  (\ref{mu(n,k)}) where $n = k +2$, and let $\displaystyle \left\{ f_k \right\}_{k=0}^\infty$ be the Fibonacci sequence (\ref{Fib}). Then for every integer $k = 0,1,2,3, \ldots$ the equality $$\mu (k+2,k) =f_k$$ is true.
\end{theorem}

Proof. When $k=0$ the assertion follows from (\ref{Fib}) and (\ref{keq0}). When $k=1$ the assertion follows from (\ref{Fib}) and (\ref{keq1}). When $k=2$ there are two canonical elements in $\mathcal{T}_4$ and these are $\mathbf{x}_1 =\langle 3,3,12,12\rangle$ and $\mathbf{x}_2 =\langle 3,5,10,12 \rangle$ (see (\ref{k2}) and \cite{umb2013}). Therefore
$$\mu (2,0) =f_0 ,\quad \mu (3,1) =f_1 \quad \textrm{and} \quad \mu (4,2) =f_2 $$

Let $k$ be an arbitrary positive integer such that $k\ge 3$ and let $A=(\alpha_{i\; j} )\in \Lambda_{k+2}^k$, $1\le i,j\le k+2$ be a canonical matrix. Then, according to Lemma  \ref{lemm1}   $\alpha_{1\, 1}=\alpha_{1\, 2}= \alpha_{2\, 1}=0$ and $\alpha_{1\, 3} =\alpha_{1\, 4} =\ldots =\alpha_{1\, n} =\alpha_{3\, 1} =\alpha_{4\, 1} =\ldots =\alpha_{n\, 1}=1$. Therefore, the following two cases are possible:

i) $\alpha_{2\, 2} =0$, i.e., $A$ is of the form
$$A= \left(
\begin{array}{ccccc}
0 & 0 & 1 & \cdots & 1 \\
0 & 0 & 1 & \cdots & 1 \\
1 & 1 &   &    &   \\
\vdots & \vdots &   & B &   \\
1 & 1 &   &   &
\end{array}
\right)
$$

We denote by $\mathcal{M}_1$ the set of all canonical matrices of this kind. Let $A$ be an arbitrary matrix of $\mathcal{M}_1$. In $A$, we remove the first and second rows and the first and second columns. We obtain the matrix $B\in \Lambda_k^{k-2}$. It is easy to see that $ B $ is the canonical matrix.

Conversely, let $B =(\beta_{i\, j} ) \in \Lambda_k^{k-2}$ ($k\ge 3$) and let $B$ be a canonical matrix. From $B$ we obtain the matrix $A=(\alpha_{i\, j})\in\Lambda_{k+2}^k$ as follows: $\alpha_{1\, 1} =\alpha_{1\, 2} =\alpha_{2\, 1}=\alpha_{2\, 2} =0$, $\alpha_{1\, j} =\alpha_{2\, j} =1$, $3\le j\le k+2$ and $\alpha_{i\, 1} =\alpha_{i\, 2} =1$, $3\le i\le k+2$. For each  $i,j\in \{3,4,\ldots ,k+2 \}$ we assume $\alpha_{i\, j} =\beta_{i-2\, j-2}$. It is easy to see that the so obtained matrix $A$ is a canonical matrix.

Therefore, $|\mathcal{M}_1 | =\mu (k,k-2)$ for any integer $k\ge 3$.

ii)  $\alpha_{2\, 2} =1$, i.e., $A$ is of the form
$$A= \left(
\begin{array}{cccccc}
0 & 0 & 1 & 1 & \cdots & 1 \\
0 & 1 & 0 & 1 & \cdots & 1 \\
1 & 0 &   &   &    &   \\
1 & 1 &   &   &    &   \\
\vdots & \vdots &   &  &  &   \\
1 & 1 &   &   &  &
\end{array}
\right)
$$

Let $\mathcal{M}_2$ be the set of all canonical matrices of this kind and let $A=(\alpha_{i\, j} )$, $\alpha_{2\, 2}=1$ be an arbitrary matrix of $\mathcal{M}_2$. We change  $\alpha_{2\, 2}$ from 1 to 0 and remove the first row and the first column of $A$. In this way we obtain a matrix of  $\Lambda_{k+1}^{k-1}$, which is easy to see that it is canonical.

Conversely, let $B=(\beta_{i\, j} ) \in \Lambda_{k+1}^{k-1}$ and let $B$ be a canonical matrix. According to Lemma \ref{lemm1} $\beta_{1\, 1} =\beta_{1\, 2} =\beta_{2\, 1} =0$. We change  $\beta_{1\, 1}$  from 0 to 1. In $B$, we add a first row and a first column and get the matrix $A=(\alpha_{i\, j} )\in \Lambda_{k+2}^k$, such that $\alpha_{1\, 1} =\alpha_{1\, 2} =\alpha_{2\, 1} =0$, $\alpha_{1\, j} =1$ for $j=3,4,\ldots ,k+2$, $\alpha_{i\, 1} =1$ for $i=3,4,\ldots ,k+2$ and $\alpha_{s+1\, t+1} =\beta_{s\, t}$ for $s,t\in \{ 1,2,\ldots ,k+1\}$. It is easy to see that the resulting matrix $A$ is canonical and  $A\in \mathcal{M}_2$.

Therefore,  $|\mathcal{M}_2 | =\mu (k+1,k-1)$ for every integer $k\ge 3$.

If $\mathcal{M}$ is the set of all canonical matrices, $\mathcal{M} \subseteq \Lambda_{k+2}^k$, then obviously
$$\mathcal{M}_1 \cap \mathcal{M}_2 =\emptyset \quad \textrm{and} \quad \mathcal{M}_1 \cup \mathcal{M}_2 =\mathcal{M} .$$

Therefore
$$\mu (k+2,k) =|\mathcal{M} |=|\mathcal{M}_1 |+|\mathcal{M}_2 |=\mu (k,k-2)+\mu (k+1,k-1)$$
for all integers $k\ge 3$, which proves the theorem.
\hfill $\Box$

\bibliographystyle{plain}
\bibliography{fib}

\end{document}